\documentclass[10pt]{article}

\usepackage[utf8]{inputenc}
\usepackage[T2A]{fontenc}
\usepackage[english]{babel}
\usepackage{amsthm, amsfonts, amsmath}
\usepackage{graphicx}
\usepackage{epstopdf}

\makeatletter
\def\th@plain{
  \thm@notefont{}
  \itshape
}
\def\th@defenition{
  \thm@notefont{}
  \normalfont
}
\def\th@remark{
  \thm@notefont{}
  \normalfont
}

\DeclareFontFamily{U}{tipa}{}
\DeclareFontShape{U}{tipa}{m}{n}{<->tipa10}{}
\newcommand{\arc@char}{{\usefont{U}{tipa}{m}{n}\symbol{62}}}%

\newcommand{\arc}[1]{\mathpalette\arc@arc{#1}}

\newcommand{\arc@arc}[2]{%
  \sbox0{$\m@th#1#2$}%
  \vbox{
    \hbox{\resizebox{\wd0}{\height}{\arc@char}}
    \nointerlineskip
    \box0
  }%
}
\makeatother

\title{Symmetries of 3-Polytopes with Fixed Edge Lengths}
\author{E.\,A.\,Morozov\thanks{National Research University Higher School of Economics, Faculty of Mathematics\endgraf gorg.morozov@gmail.com}}
\date{}

\theoremstyle{plain}
\newtheorem{lemma}{Lemma}
\newtheorem{corollary}{Corollary}
\newtheorem{theorem}{Theorem}
\newtheorem{conjecture}{Conjecture}

\theoremstyle{remark}
\newtheorem{remark}{Remark}
\newtheorem{example}{Example}

\begin{document}

\maketitle

\begin{abstract}
We consider an interesting class of combinatorial symmetries of polytopes which we call \emph{edge-length preserving combinatorial symmetries}. These symmetries not only preserve the combinatorial structure of a polytope but also map each edge of the polytope to an edge of the same length. We prove a simple sufficient condition for a polytope to realize all edge-length preserving combinatorial symmetries by isometries of ambient space. The proof of this condition uses Cauchy's rigidity theorem in an unusual way.

\smallskip


\smallskip

\end{abstract}

\section{Introduction and motivation}
Take a convex 3-dimensional polytope $P$ and consider an isometry of Euclidean 3-space mapping $P$ to itself. This isometry induces an automorphism of the combinatorial structure of the polytope. However the converse is not always true: a given combinatorial symmetry is not necessarily realizable by an isometry of ambient space. But if a combinatorial symmetry is realizable, then it maps each edge of the polytope to an edge of the same length. Such combinatorial symmetry is called \emph{edge-length preserving} (for formal definitions see Sections~\ref{sec:defnot} and~\ref{sec:symmetries}).

Our main result (Theorem~\ref{thm:main1}) states that if each 2-dimensional face of a convex 3-dimensional polytope is cyclic, i.~e. inscribed in a circle, then the polytope realizes all its own edge-length preserving combinatorial symmetries. Moreover, we prove an analogous result for convex plane graphs with cyclic bounded faces (Theorem~\ref{thm:main1forgraphs}). Note that such graphs (which are also known as \emph{circle patterns}) appear naturally in discrete complex analysis~\cite{bib:bs}, \cite{bib:kenyon}.

The proof of the main result is surprisingly short. We consider each edge-length preserving combinatorial symmetry separately and show that the given polytope realizes the symmetry by using Cauchy's rigidity theorem. Speaking informally, we prove that the given polytope is congruent to itself with nontrivial permutation of vertices.

The idea of studying combinatorial symmetries of polytopes with prescribed metric properties is not new. In a particular case (for simplicial polytopes only) edge-length preserving combinatorial symmetries appear in~\cite[\S2]{bib:gms}. Furthermore, it is known that for any 3-polytope $P$ there exists a 3-polytope $Q$ such that $P$ and $Q$ are combinatorially equivalent and $Q$ realizes all its own combinatorial symmetries (see~\cite[p. 279]{bib:mani} and~\cite[Theorem 2 in introduction]{bib:springborn}; $Q$ is called a \emph{canonical form} of the polytope $P$). This result is a simple corollary of a similar assertion about circle packings (for which M\"obius transformations are considered instead of isometries). Circle packings are related to discrete complex analysis as well. We suggest the existence of an edge-length preserving analogue of a canonical form as a conjecture (Conjecture~\ref{con:main}).

\section{Definitions and notation}\label{sec:defnot}
First recall some basic definitions concerning convex polytopes. The majority of these definitions is taken from \cite{bib:alexandrov} and \cite{bib:ziegler}.

A \emph{convex polytope} is a bounded intersection of finitely many closed half-spaces in $\mathbb R^n$. Suppose that $P$ is a convex polytope and $W$ is a subspace of $\mathbb R^n$ such that $P\subset W$. The dimension of the minimal possible such $W$ is called the \emph{dimension} of $P$ (in this case, $P$ is \emph{d-polytope}). Any hyperplane $H$ divides $\mathbb R^n$ into two closed half-spaces $H_+$ and $H_-$. The hyperplane $H$ is called \emph{supporting} if $P$ is contained in one of the half-spaces $H_+$ or $H_-$ and $H\cap P\ne\emptyset$. In this case, the intersection $H\cap P$ is a \emph{face} of the polytope $P$. Note that faces of a convex polytope are also convex polytopes. The 0-dimensional and 1-dimensional faces are called \emph{vertices} and \emph{edges} respectively. Denote by $V(P)$ and $E(P)$ the sets of vertices and edges of the convex polytope $P$ respectively and by $|V(P)|$ the number of vertices.

Define a \emph{labeled polytope} as a polytope $P$ with vertices labeled in a one-to-one fashion by the elements of some $|V(P)|$-element set. Usually this set is just $\{1,\dots,|V(P)|\}$. Note that the labeling of each face of $P$ is induced by the labeling of $P$. In what follows by a polytope we always mean a convex labeled polytope.

Two labeled polytopes are called \emph{congruent} if they are labeled by the elements of the same set and there exists an isometry of ambient space $\mathbb R^n$ which maps each vertex of the first polytope to the vertex of the second polytope with the same label.

Denote by $L(P)$ the set of all faces of a polytope $P$. The set $L(P)$, ordered by inclusion, is called the \emph{face lattice} of the polytope $P$. Two polytopes $P$ and $P'$ are called \emph{combinatorially equivalent} if there exists an order-preserving bijection $\sigma\colon L(P)\to L(P')$ such that for each vertex $v\in V(P)$ the labels of $v$ and $\sigma(v)$ are the same (in particular, the vertices of $P$ and $P'$ are labeled by the elements of the same set). In the latter case, faces $F\in L(P)$ and $\sigma(F)\in L(P')$ are called \emph{corresponding}. Note that $\sigma$ induces a combinatorial equivalence of corresponding faces.

\begin{figure}[h]
\begin{minipage}{0.23\linewidth}
\center{\includegraphics[width=1\linewidth]{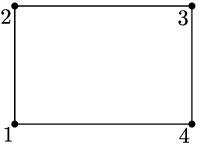}} \\ $a$
\end{minipage}
\hfill
\begin{minipage}{0.23\linewidth}
\center{\includegraphics[width=1\linewidth]{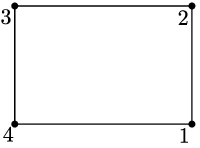}} \\ $b$
\end{minipage}
\hfill
\begin{minipage}{0.23\linewidth}
\center{\includegraphics[width=1\linewidth]{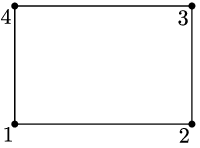}} \\ $c$
\end{minipage}
\hfill
\begin{minipage}{0.23\linewidth}
\center{\includegraphics[width=1\linewidth]{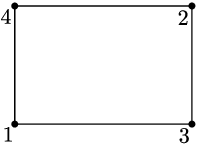}} \\ $d$
\end{minipage}
\caption{To Example~\ref{exm:4rectangles}: the definitions of congruence, combinatorial equivalence, and combinatorial symmetry of polytopes depend on the labeling of vertices.}
\label{fig:rectangles}
\end{figure}

\begin{example}[four rectangles]\label{exm:4rectangles}
Our definitions of congruence and combinatorial equivalence of polytopes depend on the labeling of vertices. For example, in Fig.~\ref{fig:rectangles} the rectangles $a$ and $b$ are congruent, the rectangles $a$ and $c$ are combinatorially equivalent but not congruent, and the rectangles $a$ and $d$ are even not combinatorially equivalent. However if we do not label the vertices, then all four rectangles become congruent.
\end{example}

The key to our proof of the main theorem is one of the most celebrated assertions about the rigidity of convex 3-polytopes.

\begin{theorem}[Cauchy, {\cite{bib:cauchy}}]\label{thm:cauchy}
If each two corresponding 2-dimensional faces of two given combinatorial equivalent labeled convex 3-polytopes are congruent, then these polytopes themselves are congruent.
\end{theorem}

Theorem~\ref{thm:cauchy} was proved by Cauchy in 1813 by the famous ``Cauchy's signs method''. This proof is widely known (see~\cite[Chapter~14]{bib:thebook} or~\cite[Theorem~24.1]{bib:omnibus}).

\begin{remark}
Actually Theorem~\ref{thm:cauchy} is usually stated in terms of polytopes without labeling of vertices. However the proof of our version of the theorem is almost literally the same.
\end{remark}

\section{Statements and proofs}\label{sec:symmetries}
Suppose that $P$ is a labeled convex polytope and $\sigma$ is an arbitrary permutation of the vertices of $P$. Let us relabel the vertices of $P$ as follows: the new label of a vertex $v$ is the former label of the vertex $\sigma^{-1}(v)$. Denote by $P_\sigma$ the obtained polytope. The permutation $\sigma$ is called a \emph{combinatorial symmetry} of $P$, if $P$ and $P_\sigma$ are combinatorially equivalent; in other words, if the permutation $\sigma\colon V(P)\to V(P)$ extends to an isomorphism $\tilde\sigma\colon L(P)\to L(P)$. In the sequel, we sometimes write just ``symmetry'' instead of ``combinatorial symmetry''. We write permutations of vertices as permutations of their labels in our examples.

The combinatorial symmetry $\sigma$ is called \emph{edge-length preserving} if for each edge $e$ of the polytope $P$ the lengths of the edges $e$ and $\tilde\sigma(e)$ are equal. We say that the polytope $P$ \emph{realizes} its combinatorial symmetry $\sigma$ if the polytopes $P$ and $P_\sigma$ are congruent, i. e., $\sigma$ extends to an isometry of $\mathbb R^3$. Furthermore, we say that the polytope $P$ \emph{realizes all its own combinatorial symmetries} if for any combinatorial symmetry $\sigma$ of $P$ the polytope $P$ realizes $\sigma$. Similarly, we say that the polytope $P$ \emph{realizes all its own edge-length preserving combinatorial symmetries} if for any edge-length preserving combinatorial symmetry $\sigma$ of $P$ the polytope $P$ realizes $\sigma$.

\begin{remark}
Obviously, the last two notions do not depend on vertices labeling.
\end{remark}

\begin{example}
Let $P$ be the rectangle in Fig.~\ref{fig:rectangles}$a$, $\sigma_1=\bigl(\begin{smallmatrix}1 & 2 & 3 & 4\\4 & 3 & 2 & 1\end{smallmatrix}\bigr)$, $\sigma_2=\bigl(\begin{smallmatrix}1 & 2 & 3 & 4\\1 & 4 & 3 & 2\end{smallmatrix}\bigr)$. Then the rectangles $P_{\sigma_1}$ and $P_{\sigma_2}$ are shown in Fig.~\ref{fig:rectangles}$b$ and Fig.~\ref{fig:rectangles}$c$ respectively. By Example~\ref{exm:4rectangles} it follows that $\sigma_1$ and $\sigma_2$ are combinatorial symmetries of $P$. Moreover, the symmetry $\sigma_1$ is edge-length preserving and $P$ realizes this symmetry.
\end{example}

\begin{example}[parallelepipeds]
A generic parallelepiped does not realize all its own edge-length preserving combinatorial symmetries. Only a rectangular parallelepiped has this property. However a generic rectangular parallelepiped does not realize all its own combinatorial symmetries. Only a cube has the latter property.
\end{example}

Recall that a polygon is called \emph{cyclic} if all its vertices lie on a circle. The main result of this paper is the following theorem.

\begin{theorem}[main theorem]\label{thm:main1}
If each 2-dimensional face of a convex 3-polytope $P$ is cyclic, then $P$ realizes all its own edge-length preserving combinatorial symmetries.
\end{theorem}

\begin{corollary}
If convex 3-polytope $P$ is inscribed in a sphere, then $P$ realizes all its own edge-length preserving combinatorial symmetries.
\end{corollary}

\begin{corollary}
If each 2-dimensional face of a convex 3-polytope $P$ is a triangle, then $P$ realizes all its own edge-length preserving combinatorial symmetries.
\end{corollary}

A finite plane graph $G$ is called \emph{convex} if each bounded face of $G$ is a convex polygon and the unbounded face of $G$ is a complement to a convex polygon. Such terms as vertices labeling, congruence, face lattice, combinatorial symmetry, edge-length preserving combinatorial symmetry, and others for convex plane graphs are defined in the same way as for 3-polytopes; one should just replace 0-dimensional, 1-dimensional, and 2-dimensional faces of 3-polytope by vertices, edges, and faces of $G$ (including the unbounded one) respectively.

\begin{example}
Like a convex polytope, a convex plane graph does not necessarily realize all its own edge-length preserving combinatorial symmetries. For example, consider a parallelogram distinct from a rectangle.
\end{example}

For convex plane graphs, the following analogue of the main theorem holds.

\begin{theorem}\label{thm:main1forgraphs}
If each bounded 2-dimensional face of a convex plane graph $G$ is cyclic, then $G$ realizes all its own edge-length preserving combinatorial symmetries.
\end{theorem}

Now proceed to the proofs. In order to prove Theorem~\ref{thm:main1}, we need some lemmas.

\begin{lemma}[Cauchy's arm lemma]\label{lem:arm}
Suppose $P=A_1A_2\dots A_n$ and $P'=A_1'A_2'\dots A_n'$ are two convex polygons such that $A_iA_{i+1}=A_i'A_{i+1}'$ for each $i=1,\dots,n-1$ and $\angle A_{j-1}A_jA_{j+1}\le\angle A_{j-1}'A_j'A_{j+1}'$ for each $j=2,\dots,n-1$. Then $A_1A_n\le A_1'A_n'$ with equality if and only if $\angle A_{j-1}A_jA_{j+1}=\angle A_{j-1}'A_j'A_{j+1}'$ for each $j=2,\dots,n-1$.
\end{lemma}

For an elementary proof of Lemma~\ref{lem:arm} see, for example,~\cite[Chapter~14]{bib:thebook}.

\begin{lemma}\label{lem:circle}
If each two corresponding sides of two cyclic convex polygons are congruent, then these polygons themselves are congruent.
\end{lemma}

The latter lemma is well-known. The article~\cite[Theorem~1.1]{bib:kss} contains a proof based on a variational principle (some historical remarks also can be found there); however, that proof is quite complicated. We give an elementary proof.

\smallskip

\noindent\textbf{Proof of Lemma~\ref{lem:circle}.} Suppose that $P=A_1A_2\dots A_n$ and $P'=A_1'A_2'\dots A_n'$ are the given polygons, $A_iA_{i+1}=A_i'A_{i+1}'$ for each $i=1,\dots,n$ (hereafter $A_{n+1}:=A_1$), $P$ is inscribed in a circle of center $O$ and radius $r$, and $P'$ is inscribed in a circle of center $O'$ and radius $r'$. We may assume that $r\le r'$ and that $A_1A_n$ is (one of) the longest side(s) of $P$. Then $O$ lies inside the angle $\angle A_{i-1}A_{i}A_{i+1}$ for each $i=2,\dots,n-1$. Indeed, otherwise the angle measure of one of the arcs $\arc{A_{i-1}A_{i+1}A_i}$ or $\arc{A_iA_{i-1}A_{i+1}}$ is at most $\pi$. Since $A_1$ and $A_n$ lie on both arcs, we have either $A_{i-1}A_i>A_1A_n$ or $A_iA_{i+1}>A_1A_n$. This contradicts the assumption made about $A_1A_n$. Thus $O$ lies inside the angle $\angle A_{i-1}A_{i}A_{i+1}$ for each $i=2,\dots,n-1$ and analogously $O'$ lies inside the angle $\angle A_{i-1}'A_i'A_{i+1}'$ for each $i=2,\dots,n-1$.

Now fix some $2\le i\le n-1$. Let us show that $\angle A_{i-1}A_iA_{i+1}\le\angle A_{i-1}'A_i'A_{i+1}'$. Consider isosceles triangles $\bigtriangleup A_{i-1}OA_i$ and $\bigtriangleup A_{i-1}'O'A_i'$. Since $r\le r'$ and $A_{i-1}A_i=A_{i-1}'A_i'$, we obtain $\angle A_{i-1}A_iO\le\angle A_{i-1}'A_i'O'$. Similarly, $\angle A_{i+1}A_iO\le\angle A_{i+1}'A_i'O'$. Since $O$ and $O'$ lie inside the angles $\angle A_{i-1}A_iA_{i+1}$ and $\angle A_{i-1}'A_i'A_{i+1}'$ respectively, we have
$$
\angle A_{i-1}A_iA_{i+1}=
\angle A_{i-1}A_iO+\angle A_{i+1}A_iO\le
\angle A_{i-1}'A_i'O'+\angle A_{i+1}'A_i'O'=
\angle A_{i-1}'A_i'A_{i+1}'.
$$
Hence $\angle A_{i-1}A_iA_{i+1}\le\angle A_{i-1}'A_i'A_{i+1}'$ for each $i=2,\dots,n-1$. Since $A_1A_n=A_1'A_n'$, from Lemma~\ref{lem:arm} it follows that in fact $\angle A_{i-1}A_iA_{i+1}=\angle A_{i-1}'A_i'A_{i+1}'$ for each $i=2,\dots,n-1$, i.~e., $P$ and $P'$ are congruent.\qed

\smallskip

\noindent\textbf{Proof of Theorem~\ref{thm:main1}.} First recall that $P$ is a \emph{labeled} polytope. Suppose $\sigma$ is an arbitrary edge-length preserving combinatorial symmetry of the polytope $P$ and $\tilde\sigma$ is an isomorphism extending $\sigma$. Then $P_\sigma$ is a labeled polytope as well. This polytope is different from $P$ in the vertices labels only. To prove the theorem it suffices to show that the polytopes $P$ and $P_\sigma$ are congruent.

Take an arbitrary 2-dimensional face $F$ of the polytope $P$ and the corresponding face $\tilde\sigma(F)$ of the polytope $P_\sigma$. Since the symmetry $\sigma$ is edge-length preserving, it follows that the corresponding sides of the polygons $F$ and $\tilde\sigma(F)$ are equal. Moreover, the polygons $F$ and $\tilde\sigma(F)$ are cyclic. It follows that $F$ and $\tilde\sigma(F)$ satisfy the conditions of Lemma~\ref{lem:circle} and therefore are congruent. Thus each two corresponding 2-dimensional faces of the polytopes $P$ and $P_\sigma$ are congruent and from Theorem~\ref{thm:cauchy} it follows that $P$ and $P_\sigma$ are congruent.\qed

\smallskip

To prove Theorem~\ref{thm:main1forgraphs} we need the following analogue of Cauchy's theorem for convex plane graphs.

\begin{lemma}[degenerate Cauchy's theorem]\label{lem:cauchyforgraphs}
If each two corresponding bounded faces of two given convex plane graphs are congruent, then the graphs themselves are congruent.
\end{lemma}

\noindent\textbf{Proof of Lemma~\ref{lem:cauchyforgraphs}.} Let $G$ and $G'$ be the given graphs. For a face (respectively, edge) $x$ of the graph $G$ denote by $x'$ the corresponding face (respectively, edge) of the graph $G'$.

Let $A$ be an arbitrary face of $G$. Consider an isometry $\tau$ such that $\tau(A)=A'$. It suffices to show that for any face $B$ of the graph $G$ we have $\tau(B)=B'$. Let $A_1,\dots,A_n$ be the faces of $G$ such that $A_1=A$, $A_n=B$, and $A_i$, $A_{i+1}$ have a common edge $e_i$ for each $i=1,\dots,n-1$. Let us prove that $\tau(A_i)=A_i'$ for $i=1,\dots,n$ by induction over $i$. For the case $i=1$ there is nothing to prove. If $\tau(A_i)=A_i'$ for some $i<n$ (as labeled polygons), then we have $\tau(e_i)=e_i'$. Denote by $A_{i+1}''$ the face $A_{i+1}'$ reflected across the line spanned by $e_i'$. Then since $\tau(e_i)=e_i'$ we have either $\tau(A_{i+1})=A_{i+1}'$ or $\tau(A_{i+1})=A_{i+1}''$. But the second case is impossible because then $\tau(A_i)=A_i'$ and $\tau(A_{i+1})=A_{i+1}''$ border upon the same side of $e_i'$, whereas $A_i$ and $A_{i+1}$ border upon the opposite sides of $e_i$. So $\tau(A_{i+1})=A_{i+1}'$ and this completes the proof.\qed

\noindent\textbf{The proof of Theorem~\ref{thm:main1forgraphs}} is almost literally the same as the proof of Theorem~\ref{thm:main1}, one should just replace the word ``polytope'' by ``convex plane graph'' and refer to Lemma~\ref{lem:cauchyforgraphs} instead of Theorem~\ref{thm:cauchy}.\qed

\section{Open problems}

Some generalizations of Theorem~\ref{thm:main1} to higher dimensions are possible but one has to generalize Cauchy's theorem to higher dimensions for the proofs. It is quite strange that we could not find such generalization in literature (although the reader can find a statement and a sketch of a proof in \cite[\S3.6.5]{bib:alexandrov} or~\cite[Theorem~27.2]{bib:pak}).

\begin{conjecture}
If each 2-dimensional face of a convex $d$-polytope $P$ is cyclic (where $d\ge 2$), then $P$ realizes all its own edge-length preserving combinatorial symmetries.
\end{conjecture}

Theorem~\ref{thm:main1} gives only a sufficient condition for a 3-polytope to realize all its own edge-length preserving combinatorial symmetries. This condition is not necessary as the following example shows.

\begin{figure}[h]
\begin{minipage}{0.3\linewidth}
\center{\includegraphics[width=1\linewidth]{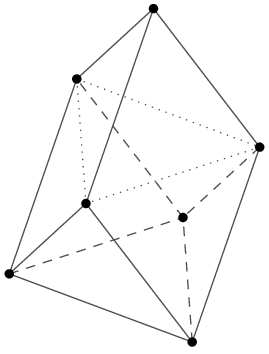}}
\end{minipage}
\hfill
\begin{minipage}{0.3\linewidth}
\center{\includegraphics[width=1\linewidth]{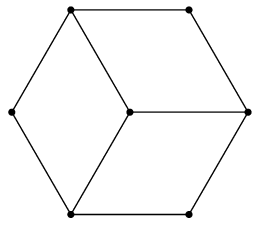}}
\end{minipage}
\caption{To Example~\ref{exm:not-nec}: a 3-polytope and a convex plane graph showing that the conditions of Theorem~\ref{thm:main1} and Theorem~\ref{thm:main1forgraphs} are not necessary ones.}
\label{fig:contra}
\end{figure}

\begin{example}\label{exm:not-nec}
Consider the polytope $P$ obtained as a union of an octahedron and a regular tetrahedron with a common face (see Fig.~\ref{fig:contra} to the left). The polytope has 4 triangular faces and 3 rhombic faces with acute angle $\pi/3$. It is easy to see that $P$ realizes all its own edge-length preserving combinatorial symmetries. Nevertheless not all faces of $P$ are cyclic. To obtain an analogous example for convex plane graphs, take a regular hexagon divided into three rhombi again with acute angle $\pi/3$ (see Fig.~\ref{fig:contra} to the right).
\end{example}

Taking into account these examples, the following conjecture seems to be very interesting.

\begin{conjecture}\label{con:main}
For each convex 3-polytope $P$ there exists a convex 3-polytope $Q$ such that $Q$ is combinatorially equivalent to $P$, the lengths of each two corresponding edges of $P$ and $Q$ are equal, and $Q$ realizes all its own edge-length preserving combinatorial symmetries.
\end{conjecture}

Note that the latter conjecture is not true for convex plane graphs!

\begin{figure}[h]
\begin{minipage}{0.3\linewidth}
\center{\includegraphics[width=1\linewidth]{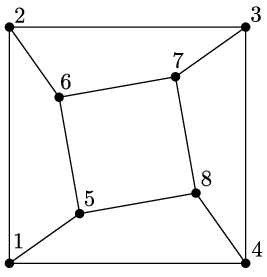}}
\end{minipage}
\hfill
\begin{minipage}{0.4\linewidth}
\center{\includegraphics[width=1\linewidth]{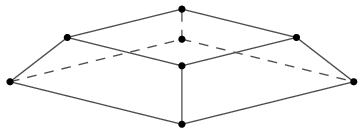}}
\end{minipage}
\caption{To Example~\ref{exm:cube}: a counterexample to a direct analogue of Conjecture~\ref{con:main} for convex plane graphs and a 3-polytope fixing this counterexample.}
\label{fig:cube}
\end{figure}

\begin{example}\label{exm:cube}
Consider a convex plane graph $G$ in Fig.~\ref{fig:cube} to the left. This graph consists of two squares $(1234)$ and $(5678)$ with a common center. Clearly, if the angle between the sides $(12)$ and $(56)$ of the squares is positive but small enough, then all bounded faces of $G$ are convex quadrilaterals, i. e., $G$ is indeed a convex plane graph.

Let us prove that there is no convex plane graph $G'$ such that $G'$ is combinatorially equivalent to $G$, the lengths of each two corresponding edges of $G$ and $G'$ are equal, and $G'$ realizes all its own edge-length preserving combinatorial symmetries. Assume the converse. Then $G'$ realizes the symmetry $\bigl(\begin{smallmatrix}1 & 2 & 3 & 4 & 5 & 6 & 7 & 8\\2 & 1 & 4 & 3 & 6 & 5 & 8 & 7\end{smallmatrix}\bigr)$. It follows that $\angle(215)=\angle(126)$ in $G'$. Similarly, $\angle(215)=\angle(326)=\angle(437)=\angle(148)=\angle(126)=\angle(237)=\angle(348)=\angle(415)=\pi/4$ in $G'$ (the last equality holds since the sum of all eight angles is $2\pi$). Hence the quadrilaterals $(1234)$ and $(5678)$ in $G'$ are squares with a common center and parallel sides. But if the sizes of these squares are fixed, then the length of the edge $(15)$ grows as the angle between the lines $(12)$ and $(56)$ grows. Thus the edge $(15)$ cannot have the same length in both graphs $G$ and $G'$.

However the graph $G$ can be considered as a ``degenerate'' polytope with all the faces lying in one plane. Under this point of view, there exists a \emph{3-polytope} $P$ such that $P$ has the same face lattice as $G$, the lengths of the corresponding edges of $G$ and $P$ are equal, and $P$ realizes all its own edge-length preserving combinatorial symmetries. This polytope is square frustum (see Fig.~\ref{fig:cube} to the right).
\end{example}

\section*{Acknowledgments}
I am grateful to my scientific supervisor M.\,B.\,Skopenkov for constant attention to this work and great help in preparing this text. I am also grateful to I.\,Kh.\,Sabitov and M.\,S.\,Tyomkin for interesting discussions.

\end{document}